\documentclass[a4paper,12pt]{amsart}
\usepackage{amsfonts}
\usepackage{amssymb}
\usepackage{ifthen}
\usepackage{graphicx}
\usepackage{color}
\nonstopmode \numberwithin{equation}{section}
\newtheorem{thm}{Theorem}
\newtheorem{lem}{Lemma}
\newtheorem{cor}{Corollary}[section]

\newtheorem{cl}{Claim}
\newtheorem{ca}{Case}
\newtheorem{sca}{Subcase}
\newtheorem{scl}{Subclaim}
\newtheorem{conj}{Conjecture}

\theoremstyle{definition}
\newtheorem{defn}{Definition}

\newtheorem{op}[equation]{Open Problem}
\newtheorem{ques}[equation]{Question}
\newtheorem{rem}{Remark}[section]
\newtheorem{exam}[equation]{Example}

\newcounter {own}
\def\theown {\thesection       .\arabic{own}}

\newenvironment{pf}[1][]{%
 \vskip 3mm
 \noindent
 \ifthenelse{\equal{#1}{}}%
  {{\slshape Proof. }}%
  {{\slshape #1.} }%
 }%
{\qed\bigskip}

\newcounter{alphabet}
\newcounter{tmp}
\newenvironment{Thm}[1][]{\refstepcounter{alphabet}%
\bigskip%
\noindent%
{\bf Theorem \Alph{alphabet}}%
\ifthenelse{\equal{#1}{}}{}{ (#1)}%
{\bf .} \itshape}{\vskip 8pt}

\makeatletter
\newcommand{\Ref}[1]{\@ifundefined{r@#1}{}{\setcounter{tmp}{\ref{#1}}\Alph{tmp}}}
\makeatother

\newenvironment{Lem}[1][]{\refstepcounter{alphabet}%
\bigskip%
\noindent%
{\bf Lemma \Alph{alphabet}}%
{\bf .} \itshape}{\vskip 8pt}

\newcommand{\ID}{{\mathbb D}}

\def\be{\begin{equation}}
\def\ee{\end{equation}}

\newcommand{\bee}{\begin{enumerate}}
\newcommand{\eee}{\end{enumerate}}

\newcommand{\blem}{\begin{lem}}
\newcommand{\elem}{\end{lem}}
\newcommand{\bthm}{\begin{thm}}
\newcommand{\ethm}{\end{thm}}
\newcommand{\bcor}{\begin{cor}}
\newcommand{\ecor}{\end{cor}}
\newcommand{\beg}{\begin{exam}}
\newcommand{\eeg}{\end{exam}}
\newcommand{\begs}{\begin{examples}}
\newcommand{\eegs}{\end{examples}}
\newcommand{\bdefe}{\begin{defn}}
\newcommand{\edefe}{\end{defn}}
\newcommand{\bprob}{\begin{prob}}
\newcommand{\eprob}{\end{prob}}
\newcommand{\bques}{\begin{ques}}
\newcommand{\eques}{\end{ques}}
\newcommand{\bei}{\begin{itemize}}
\newcommand{\eei}{\end{itemize}}
\newcommand{\bcon}{\begin{conj}}
\newcommand{\econ}{\end{conj}}
\newcommand{\bop}{\begin{op}}
\newcommand{\eop}{\end{op}}

\newcommand{\bca}{\begin{ca}}
\newcommand{\eca}{\end{ca}}
\newcommand{\bsca}{\begin{sca}}
\newcommand{\esca}{\end{sca}}

\newcommand{\bcl}{\begin{cl}}
\newcommand{\ecl}{\end{cl}}

\newcommand{\bscl}{\begin{scl}}
\newcommand{\escl}{\end{scl}}

\newcommand{\bcons}{\begin{conjs}}
\newcommand{\econs}{\end{conjs}}
\newcommand{\bprop}{\begin{propo}}
\newcommand{\eprop}{\end{propo}}
\newcommand{\br}{\begin{rem}}
\newcommand{\er}{\end{rem}}
\newcommand{\brs}{\begin{rems}}
\newcommand{\ers}{\end{rems}}
\newcommand{\bo}{\begin{obser}}
\newcommand{\eo}{\end{obser}}
\newcommand{\bos}{\begin{obsers}}
\newcommand{\eos}{\end{obsers}}
\newcommand{\bpf}{\begin{pf}}
\newcommand{\epf}{\end{pf}}
\newcommand{\ba}{\begin{array}}
\newcommand{\ea}{\end{array}}
\newcommand{\beq}{\begin{eqnarray}}
\newcommand{\beqq}{\begin{eqnarray*}}
\newcommand{\eeq}{\end{eqnarray}}
\newcommand{\eeqq}{\end{eqnarray*}}

\newcounter{minutes}
\divide\time by 60
\newcounter{hours}
\multiply\time by 60 \addtocounter{minutes}{-\time}

\begin{document}

\bibliographystyle{amsplain}
\title{Linear connectivity, Schwarz-Pick lemma and univalency criteria for planar harmonic mappings
}

\def\thefootnote{}
\footnotetext{ \texttt{\tiny File:~\jobname .tex,
          printed: \number\day-\number\month-\number\year,
          \thehours.\ifnum\theminutes<10{0}\fi\theminutes}
} \makeatletter\def\thefootnote{\@arabic\c@footnote}\makeatother

\author{Sh. Chen }
\address{Sh. Chen, Department of Mathematics and Computational
Science, Hengyang Normal University, Hengyang, Hunan 421008,
People's Republic of China.} \email{mathechen@126.com}

\author{S. Ponnusamy $^\dagger $
}
\address{S. Ponnusamy,
Indian Statistical Institute (ISI), Chennai Centre, SETS (Society
for Electronic Transactions and security), MGR Knowledge City, CIT
Campus, Taramani, Chennai 600 113, India. }
\email{samy@isichennai.res.in, samy@iitm.ac.in}

\author{A. Rasila }
\address{A. Rasila, Department of Mathematics and Systems Analysis, Aalto University, P. O. Box 11100, FI-00076 Aalto,
 Finland.} \email{antti.rasila@iki.fi}

\author{X. Wang  $^\dagger {}^\dagger$
}

\address{X. Wang, Department of Mathematics,
Hunan Normal University, Changsha, Hunan 410081, People's Republic
of China.} \email{xtwang@hunnu.edu.cn}

\keywords{Harmonic  mapping, linearly connected domain, $\alpha$-close-to-convex function, John constant, univalency.\\
$^\dagger$ {\tt This author is on leave from the
Department of Mathematics,
Indian Institute of Technology Madras, Chennai-600 036, India.\\
${}^\dagger {}^\dagger$ {\tt Corresponding author.} } }


\begin{abstract} In this paper, we first establish the Schwarz-Pick lemma
of higher-order  and  apply it to obtain a univalency criteria for
planar harmonic mappings. Then  we discuss distortion theorems,
Lipschitz continuity  and univalency of planar harmonic mappings defined in the unit disk with
linearly connected images.
\end{abstract}


\maketitle
\pagestyle{myheadings}
\markboth{ Sh. Chen, S. Ponnusamy, A. Rasila, X. Wang }{Linear connectivity and Schwarz-Pick lemma for planar harmonic mappings}

\section{Introduction and main results }\label{csw-sec1}
Let $f$ be a complex-valued function defined on a simply connected subdomain $D$ of the complex plane
$\mathbb{C}$ and $f=h+\overline{g}$ its decomposition (unique up to an additive constant),
where $h$ and $g$ are analytic in $D$. It is convenient to choose the additive constant in such
a way that $g(0)=0$. In this case, the decomposition is unique, and it is called the canonical
decomposition. Since the {\it Jacobian} $J_f$ of $f$ is given by
$$J_f=|f_{z}|^2-|f_{\overline{z}}|^2 =|h'|^2-|g'|^2,
$$
$f$ is locally univalent and sense-preserving in $D$ if and only if
$|g'(z)|<|h'(z)|$ in $D$. The
(second) complex dilatation $\omega =g'/h'$ of a sense-preserving harmonic mapping $f$ has the
property that $|\omega(z)|<1$ in $D$ (see \cite{Lewy}).
We refer to \cite{Clunie-Small-84,Du,PonRasi2013} for basic results in the theory of planar
harmonic mappings.

We first recall that the classical Schwarz Lemma for analytic
functions $f$ of $\mathbb{D}$ into itself as follows:
\be\label{eq-CPR1}
|f'(z)|\leq\frac{1-|f(z)|^{2}}{1-|z|^{2}}.
\ee
In 1920, Sz\'asz \cite{S} extended the inequality (\ref{eq-CPR1}) to
the following estimate involving higher order derivatives:

\be\label{eq-CPR2}
|f^{(2m+1)}(z)|\leq\frac{(2m+1)!}{(1-|z|^{2})^{2m+1}}\sum_{k=0}^{m}{m\choose k}^{2}|z|^{2k},
\ee
where $m\in\{1,2,\ldots\}.$ In 1985, Ruscheweyh \cite{AW, R} improved (\ref{eq-CPR2}) to the following form:
\be\label{eq-CPR3}
|f^{(n)}(z)|\leq\frac{n!(1-|f(z)|^{2})}{(1-|z|)^{n}(1+|z|)}.
\ee

In 1989, Colonna established an analogue of the Schwarz-Pick lemma
for planar harmonic mappings.

\begin{Thm}{\rm (\cite[Theorem 3]{Co})}\label{Co}
Let $f$ be a harmonic mapping of $\mathbb{D}$ into $\mathbb{C}$ such
that $\sup_{z\in\mathbb{D}}|f(z)|\leq M$, where $M$ is a positive
constant. Then for $z\in\mathbb{D}$,
$$\Lambda_{f}(z)\leq\frac{4M}{\pi}\frac{1}{1-|z|^{2}}.
$$
This estimate is sharp and all the extremal functions are
$$f(z)=\frac{2M\alpha }{\pi}\arg \left ( \frac{1+\psi(z)}{1-\psi(z)}\right),
$$
where $|\alpha|=1$ and $\psi$ is a conformal automorphism of $\mathbb{D}$.
\end{Thm}

Analogously to the inequality (\ref{eq-CPR3}), Chen, Ponnusamy and Wang \cite{CPW1} established the
higher order derivatives of harmonic mappings as follows.

\begin{Thm}{\rm (\cite[Corollary 3.1]{CPW1})}\label{ThmA}
Let $f$ be a harmonic mapping of $\mathbb{D}$ into $\mathbb{C}$ such
that $\sup_{z\in\mathbb{D}}|f(z)|\leq M$, where $M$ is a positive
constant. Then for $n\geq1$ and $z\in\mathbb{D}$,
$$\left|\frac{\partial^{n}f}{\partial z^{n}}(z) \right|\leq\frac{n!M}{(1-|z|)^{n+1}}
~\mbox{and}~\left|\frac{\partial^{n}f }{\partial
\overline{z}^{n}}(z)\right|\leq\frac{n!M}{(1-|z|)^{n+1}}.
$$
\end{Thm}

The following result is a generalization of Theorem \Ref{Co} and an improvement of Theorem \Ref{ThmA}.

\begin{thm}\label{thmx}
Let  $f$ be a harmonic mapping of $\mathbb{D}$ into $\mathbb{C}$
such that $\sup_{z\in\mathbb{D}}|f(z)|\leq M$, where $M$ is a
positive constant. Then  for $n\geq1$,
\beq\label{eq-thmx}
\left|\frac{\partial^{n} f}{\partial z^{n}}(z)\right|
+\left|\frac{\partial^{n} f}{\partial \overline{z}^{n}}(z)\right|\leq\frac{n!4M}{\pi}\frac{1}{(1-|z|)^{n}(1+|z|)},
\eeq
where  $z\in\mathbb{D}$. The estimate of {\rm (\ref{eq-thmx})} is
sharp at $z=0$.
\end{thm}

From Theorem \ref{thmx}, we get the following  result as well.


\begin{cor}\label{Cor-1}
Let $f$ be a analytic function in $\mathbb{D}$. Then
$$|f^{(n)}(z)|\leq\frac{n!4\sup_{\zeta\in\mathbb{D}}| {\rm Re}f(\zeta)|}{\pi}
\frac{1}{(1-|z|)^{n}(1+|z|)}.
$$
\end{cor}

Let $\mathcal{H}$ denote all non-constant harmonic
mappings in $\mathbb{D}$. We use $\mathcal{S}_{HU}$ to
denote all the univalent harmonic mappings  in $\mathbb{D}$. For
$f\in\mathcal{H}$, let
$$M_{f}=\sup_{z\in\mathbb{D}}\Lambda_{f}(z),~m_{f}=\sup_{z\in\mathbb{D}}\lambda_{f}(z)~\mbox{ and }~\mu_{f}=\frac{M_{f}}{m_{f}},
$$
where
$$\Lambda_{f}(z)=\max_{0\leq \theta\leq 2\pi}|f_{z}(z)+e^{-2i\theta}f_{\overline{z}}(z)|
=|f_{z}(z)|+|f_{\overline{z}}(z)|
$$
and
$$\lambda_{f}(z)=\min_{0\leq \theta\leq 2\pi}|f_{z}(z)+e^{-2i\theta}f_{\overline{z}}(z)|
=\big | \, |f_{z}(z)|-|f_{\overline{z}}(z)|\, \big |.
$$
Then $J_{f}=\lambda_{f}\Lambda_{f}$ if $J_f\geq 0$. Moreover, for ${\mathcal T} =\mathcal{H}\setminus\mathcal{S}_{HU},$
we define the {\it harmonic John constant} $\gamma$ by
\be\label{johnconst}
\gamma=\inf_{f\in {\mathcal T}}\mu_{f}.
\ee
On the studies of John constant for analytic functions, see \cite{J,Y}.

\begin{thm}\label{thm-3}
Let $\gamma$ be the harmonic John constant as in \eqref{johnconst}. Then
$e^{\frac{\pi}{2}}\leq\gamma\leq e^{\pi}.$
\end{thm}

Proofs for Theorems \ref{thmx} and \ref{thm-3} will be given in Section \ref{csw-sec2}.

A domain $\Omega\subset\mathbb{C}$ is said to be {\it $M$-linearly
connected} if there exists a positive constant $M\in[1,\infty)$ such
that any two points $z, w\in \Omega$ are joined by a path
$\gamma\subset \Omega$ with
$$\ell(\gamma)\leq M|z-w|,\text{ where }\ell(\gamma)=\inf\left\{\int_{\gamma}|dz|:\, \gamma\subset \Omega\right\}.
$$
It is not difficult to verify that a $1$-linearly
connected domain is convex. For extensive discussions on linearly
connected domains, see \cite{A,CPW,CPW-3,CH,H,P}.

Let  ${\mathcal S}_{H}$ denote the class of all sense-preserving
planar harmonic univalent mappings $f=h+\overline{g}$ defined in
$\mathbb{D}$, where $h$ and $g$ are analytic function in $\ID$
normalized in a standard form: $h(0)=g(0)=0$ and $h'(0)=1$, see
\cite{Clunie-Small-84, Du}. If $g(z)$ is identically zero on the
decomposition of $f(z)$, then the class ${\mathcal S}_{H}$ in this
case reduces to the classical family $\mathcal S$ of normalized
analytic univalent functions
$h(z)=z+\sum_{n=2}^{\infty}a_{n}z^{n}$ in $\ID$. Let ${\mathcal
S}_H^{0}=\{f=h+\overline{g} \in {\mathcal S}_H: \,g'(0)=0 \} $. It
is well-known that the family ${\mathcal S}_H^{0}$ is normal and compact (see
\cite{Clunie-Small-84,Du}).


A family ${\mathcal L}$ of locally univalent harmonic mappings is called a linearly invariant family (LIF) if
for any $f=h+\overline{g}\in {\mathcal L}$ its Koebe transformation
$$K(z):=K_\phi(f)(z):=\frac{f(\phi(z)) - f(\phi(0))}{\phi'(0) h'(\phi(0))}
$$
belongs to ${\mathcal L}$ for all analytic automorphisms $\phi$ of the disk $\ID$.
The family ${\mathcal L}$ is called an affinely and linearly invariant family (ALIF) if it is LIF and
for any $f\in {\mathcal L}$ its affine transformation
$$F_\mu(z):=F_\mu(f)(z)=\frac{f(z)+\mu\overline{f(z)}}{1+\mu g'(0)}
$$
belongs to ${\mathcal L}$ for all $\mu\in \ID$. The classical order of the family ${\mathcal L}$
is defined as
$${\rm ord}\,{\mathcal L}:=\sup\{|a_2(h)|:\, f\in {\mathcal L}\}.
$$
(see \cite{Sheil-Small-90,KneStarSzy2011}). These transformations are instrumental
in the investigation of distortion theorem and Lipschitz continuity of harmonic mappings $f\in{\mathcal
S}_H^{0}.$



\begin{thm}\label{thm-2}
Let $f=h+\overline{g}\in{\mathcal S}_{H}^{0}$, where $h$ and $g$ have the form
$$h(z)=z+\sum_{n=2}^{\infty}a_{n}z^{n} ~\mbox{ and }~ g(z)=\sum_{n=2}^{\infty}b_{n}z^{n}.
$$
Then we have the following:
\bee
\item[{\rm (I)}] There is a positive constant $c_{1}<\infty$ such that for
$\xi\in\partial\mathbb{D}$ and $0\leq\rho\leq r<1$,
$$\Lambda_{f}(r\xi)\geq\frac{1}{2^{1+c_{1}}}\Lambda_{f}(\rho\xi)\left(\frac{1-r}{1-\rho}\right)^{c_{1}-1}.
$$
\item[{\rm (II)}] Furthermore, if $\Omega=f(\mathbb{D})$ is a $M$-linearly connected domain
with  $|\omega (z)|\leq c<1$, then $|b_{2}|\leq c/2$. The estimate
of $|b_{2}|$ is sharp, and the extremal function is
$f(z)=z+\frac{c}{2}\overline{z}^{2}.$

\item[{\rm (III)}] Moreover, if $\Omega=f(\mathbb{D})$ is a $M$-linearly connected
domain with $|\omega (z)|\leq c<\frac{1}{2M+1},$ then  there is a
positive constant $c_{2}$ and $c_{3}<2$ such that for
$\zeta_{1},\zeta_{2}\in\partial\mathbb{D}$,
$$|f(\zeta_{1})-f(\zeta_{2})|\geq c_{2}|\zeta_{1}-\zeta_{2}|^{c_{3}}
$$
and for $n\geq2,$
$$|a_{n}|+|b_{n}|\leq n,
$$
where $c_{2}$ depends only on $M$.
\eee
\end{thm}

We remark that Theorem \ref{thm-2}  is a generalization of
\cite[Theorem 5.7]{P}. We conjecture that $c_1$ in Theorem \ref{thm-2}(I) could be taken as
$c_1=5/2$.
Moreover, further computations suggest the following.

\begin{conj}
Suppose that  $f=h+\overline{g}\in{\mathcal S}_{H}^{0}$ and  $\Omega=f(\mathbb{D})$ is a
$M$-linearly connected domain. Then there is a positive constant
$c_{4}<2$ such that for $\xi\in\partial\mathbb{D}$ and
$0\leq\rho\leq r<1$,
$$\Lambda_{f}(r\xi)\geq\frac{1}{8}\Lambda_{f}(\rho\xi)\left(\frac{1-r}{1-\rho}\right)^{c_{4}-1}.
$$
\end{conj}

\begin{defn}\label{def-1} Let $\alpha\in[0,1)$. A univalent
analytic function $f$ is called {\it $\alpha$-close-to convex} if
there is a univalent and convex analytic function $\phi$ such that
$$\left|\arg[f'(z)/\phi'(z)]\right|\leq\frac{\alpha\pi}{2} ~\mbox{ for $z\in\mathbb{D}$}.
$$
\end{defn}

It is known that \cite[Propostion 5.8]{P} the range of every $\alpha$-close-to convex function is linearly connected.

A domain $D$ is convex in the horizontal direction (CHD) if every line parallel to the real axis
has a connected intersection with $D$. One of the beautiful results of Clunie and Sheil-Small \cite[Theorem 5.3]{Clunie-Small-84} states that
if $f=h+\overline{g}$ is a harmonic function that is locally univalent
in $\ID$ (i.e., $|\omega(z)|<1$ for all $z\in \ID$), then the function $F=h-g$ is an analytic univalent mapping
of $\ID$ onto a CHD domain if and only if $f=h+\overline{g}$ is a univalent mapping of $\ID$ onto a CHD domain.
It is easy to establish a similar result for functions which are convex in the other directions.

In \cite{CH}, the authors  discussed the relationship between linear
connectivity of the images of $\mathbb{D}$ under the planar harmonic
mappings $f=h+\overline{g}$ and under their corresponding
analytic counterparts $h$.  The following result is a generalization of the
shearing theorem of Clunie and Sheil-Small \cite{Clunie-Small-84}.

\begin{thm}\label{thm-1}
Fix $\alpha\in[0,1)$, and let $f=h+\overline{g}$ be a harmonic mapping,
where $h$ and $g$ are analytic in $\ID$.
\bee
\item[{\rm (I)}] If $h- g$ {\rm (or $h+g$)} is  $\alpha$-close-to convex and $|\omega (z) | \leq M_{1}$ for $z\in\ID$,
then $h$ is univalent and $h(\mathbb{D})$ is $M_{2}$-linearly connected domain, where
$$ M_{1}<\frac{\cos\frac{\alpha\pi}{2}}{1+\cos\frac{\alpha\pi}{2}} ~\mbox{ and }~
M_{2}=\frac{1}{\cos\frac{\alpha\pi}{2}-M_{1}(1+\cos\frac{\alpha\pi}{2})}.
$$

\item[{\rm (II)}]  If $h- g$ {\rm (or $h+g$)} is  $\alpha$-close-to convex and $|\omega (z)| \leq
M_{3} $ for $z\in\ID$, then $f_{\theta}=h+e^{i\theta}\overline{g}$ is a $K$-quasiconformal
harmonic mapping and $f_{\theta}(\mathbb{D})$ is $M_{4}$-linearly
connected domain, where $\theta\in[0,2\pi)$,
$$K=\frac{1+M_{3}}{1-M_{3}},~ M_{3}<\frac{\cos\frac{\alpha\pi}{2}}{2+\cos\frac{\alpha\pi}{2}} ~\mbox{ and }~
M_{4}=\frac{1+M_{3}}{\cos\frac{\alpha\pi}{2}-M_{3}(2+\cos\frac{\alpha\pi}{2})}.
$$
\eee
\end{thm}

The proofs of Theorems \ref{thm-2} and \ref{thm-1} will be presented
in Section \ref{csw-sec3}.

\section{The Schwarz lemma of higher-order and a univalenncy criterion for harmonic mappings }\label{csw-sec2}


\begin{Lem}{\rm (\cite[Lemma 1]{CPW0} $\mbox{or}$ \cite[Theorem 1.1]{CPW-BMMSC2011})}\label{LemA}
Let $f$ be a harmonic mapping of $\mathbb{D}$ into $\mathbb{C}$ such
that $\sup_{z\in\mathbb{D}}|f(z)|\leq M$ and
$f(z)=\sum_{n=0}^{\infty}a_{n}z^{n}+\sum_{n=1}^{\infty}\overline{b}_{n}\overline{z}^{n}$,
where $M$ is a positive constant. Then $|a_{0}|\leq M$ and for all
$n\geq 1,$
\beq\label{eq-A}
|a_{n}|+|b_{n}|\leq \frac{4M}{\pi}.
\eeq
The estimate of {\rm (\ref{eq-A})} is sharp, all the extremal
functions are $f(z)\equiv M$ and
$$f_{n}(z)=\frac{2M\alpha}{\pi}\arg\left(\frac{1+\beta z^{n}}{1-\beta z^{n}}\right),
$$
where $|\alpha|=|\beta|=1.$
\end{Lem}

\subsection*{Proof of Theorem \ref{thmx}}
 Let $f=h+\overline{g}$, where $h$ and $g$ are analytic in
$\mathbb{D}$. For any fixed $z\in\mathbb{D}$, let
$\phi(\zeta)=\frac{\zeta+z}{1+\zeta\overline{z}}$, where
$\zeta\in\mathbb{D}.$ For $\zeta\in\mathbb{D},$ set
$$F(\zeta)=f(\phi(\zeta))=h(\phi(\zeta))+\overline{g(\phi(\zeta))}=
\sum_{k=0}^{\infty}a_{k}\zeta^{k}+\sum_{k=1}^{\infty}\overline{b}_{k}\overline{\zeta}^{k}.
$$
Using Lemma \Ref{LemA} and the well-known formula (cf. \cite{AW})
$$\frac{h^{(n)}(z)(1-|z|^{2})^{n}}{n!}=\sum_{k=1}^{n}{n-1\choose n-k}\overline{z}^{n-k}a_{k}
$$
and
$$\frac{g^{(n)}(z)(1-|z|^{2})^{n}}{n!}=\sum_{k=1}^{n}{n-1\choose n-k}\overline{z}^{n-k}b_{k},
$$
we get
\beq
\nonumber
\frac{(|h^{(n)}(z)|+|g^{(n)}(z)|)(1-|z|^{2})^{n}}{n!}&\leq&\sum_{k=1}^{n}{n-1\choose n-k}|z|^{n-k}(|a_{k}|+|b_{k}|)\\
\nonumber &\leq&\frac{4M}{\pi}\sum_{k=1}^{n}{n-1\choose n-k}|z|^{n-k}\\
\nonumber &=& \frac{4M}{\pi}(1+|z|)^{n-1},
\eeq
which gives
$$|h^{(n)}(z)|+|g^{(n)}(z)|\leq\frac{n!4M}{\pi}\frac{(1+|z|)^{n-1}}{(1-|z|^{2})^{n}}.
$$
The inequality (\ref{eq-thmx}) follows. Sharpness at $z=0$ is a consequence of the sharpness part of
Lemma \Ref{LemA}. So we omit the details.
\qed

\subsection*{Proof of Corollary \ref{Cor-1}} For $z\in\mathbb{D}$, let $u(z)=\mbox{Re}(f(z))$. Without loss
of generality, we assume that $\sup_{\zeta\in\mathbb{D}}| {\rm
Re}f(\zeta)|<\infty.$ By Theorem \ref{thmx} and elementary
calculations, we get
$$|f^{(n)}(z)|=\left|\frac{\partial^{n} u}{\partial z^{n}}\right|+\left|\frac{\partial^{n} u}{\partial
\overline{z}^{n}}\right|\leq\frac{n!4\sup_{\zeta\in\mathbb{D}}| {\rm Re}f(\zeta)|}{\pi} \frac{(1+|z|)^{n-1}}{(1-|z|^{2})^{n}}.
$$
\qed

\begin{Lem}{\rm (\cite[Corollary 4.1]{Be})}\label{Lem-Be}
Let $f$ be a non-constant analytic function in $\mathbb{D}$. If
$\|f\|\leq1,$ then $f$ is univalent in $\mathbb{D}$, where
$$\|f\|=\sup_{z\in\mathbb{D}}\left\{(1-|z|^{2})\Big|\frac{f''(z)}{f'(z)}\Big|\right\}.
$$
\end{Lem}

\subsection*{Proof of Theorem \ref{thm-3}} We first prove the left part.  For $\theta\in[0,2\pi)$, let
$F_{\theta}=h+e^{i\theta}g$. For $z\in\mathbb{D}$, set
$$H_{\theta}(z)=\log\frac{F'_{\theta}(z)}{\sqrt{M_{f}m_{f}}}.
$$
Then for $z\in\mathbb{D}$,
\be\label{eq-r11}
\mbox{Re}(H_{\theta}(z))\leq\frac{1}{2}\log\frac{M_{f}}{m_{f}}.
\ee
By Corollary \ref{Cor-1}, we have

\be\label{eq-r12}
(1-|z|^{2})|H'_{\theta}(z)|
\leq\frac{4}{\pi}\sup_{z\in\mathbb{D}}\left|\mbox{Re}(H_{\theta}(z))\right|.
\ee
Applying (\ref{eq-r11}) and (\ref{eq-r12}), we obtain
$$(1-|z|^{2})|H'_{\theta}(z)|\leq\frac{2}{\pi}\log\frac{M_{f}}{m_{f}},
$$ which gives
$$(1-|z|^{2})\left|\frac{h''(z)+e^{i\theta}g''(z)}{h'(z)+e^{i\theta}g'(z)}\right|\leq\frac{2}{\pi}\log\frac{M_{f}}{m_{f}}.
$$
By using Lemma \Ref{Lem-Be}, if
$$\frac{2}{\pi}\log\frac{M_{f}}{m_{f}}\leq1,
$$
then for all $\theta\in[0,2\pi)$ the function $F_{\theta}$ is univalent, which implies that
$f=h+\overline{g}$ is univalent. Hence
$$\frac{M_{f}}{m_{f}}> e^{\frac{\pi}{2}},
$$
which yields that $\gamma\geq e^{\frac{\pi}{2}}.$

Next we come to prove the right part. Since all the analytic
functions defined in $\mathbb{D}$ are harmonic, by
\cite[Theorem]{Y}, we see that $\gamma\leq e^{\pi}.$ The proof of
the theorem is complete. \qed

\section{Linear connectivity and univalency criterion for harmonic mappings}\label{csw-sec3}

\subsection*{Proof of Theorem \ref{thm-2}}
 We first prove (I). For every constant
$\mu\in\mathbb{D}$, consider the affine mappings $F_{\mu}=f+\mu\overline{f}$, where
$f=h+\overline{g}\in{\mathcal S}_{H}^{0}$. Clearly,  $F_{\mu}\in{\mathcal S}_{H}$.  For a fixed
$\zeta\in\mathbb{D}$, we next consider the Koebe transform of $F_{\mu}$ given by
$$K(z)=\frac{F_{\mu}(\frac{z+\zeta}{1+\overline{\zeta}z})-F_{\mu}(\zeta)}{(1-|\zeta|^{2})(h'(\zeta)+\mu g'(\zeta))}
=H(z)+\overline{G(z)},
$$
which again belongs to ${\mathcal S}_{H}.$ By elementary calculations, we get
$$H(z)=z+A_{2}(\zeta)z^{2}+A_{3}(\zeta)z^{3}+\cdots,
$$
where
$$A_{2}(\zeta)=\frac{1}{2}(1-|\zeta|^{2})\frac{h''(\zeta)+\mu g''(\zeta)} {h'(\zeta)+\mu g'(\zeta)}-\overline{\zeta}.
$$
By \cite[p.~87 and p.~96]{Du}, we know that
$|A_{2}(\zeta)|$ is bounded. 
 Without loss of generality, we assume
$|A_{2}(\zeta)|\leq c_{1}<\infty,$ which implies
\beq\label{eqt-7}
\nonumber
\left|\frac{\partial}{\partial\rho}\log\big[(1-\rho^{2})(h'(\rho\xi)+\mu
g'(\rho\xi))\big]\right| &=&\left|\frac{h''(\rho\xi)+\mu
g''(\rho\xi)}{h'(\rho\xi)+\mu
g'(\rho\xi)}-\frac{2\rho\overline{\xi}} {1-\rho^{2}}\right|\\
\nonumber & \leq&\frac{2c_{1}}{1-\rho^{2}},
\eeq
where $\xi\in\partial\mathbb{D}$. Integration leads to
\beq\label{eqt-8}
\frac{(1-r^{2})|h'(r\xi)+\mu g'(r\xi)|}{(1-\rho^{2})|h'(\rho\xi)+\mu
g'(\rho\xi)|} \geq
\left(\frac{1-r}{1+r}\cdot\frac{1+\rho}{1-\rho}\right)^{c_{1}},
\eeq
which gives
\beq\label{eqt-9}
|h'(r\xi)+\mu g'(r\xi)|\geq|h'(\rho\xi)+\mu
g'(\rho\xi)|\left(\frac{1-r}{1-\rho}\right)^{c_{1}-1}
\left(\frac{1+\rho}{1+r}\right)^{c_{1}+1}.
\eeq
By (\ref{eqt-9}), we have
\beq\label{eqt-10}
\Lambda_{f}(r\xi)\geq|h'(\rho\xi)+\mu
g'(\rho\xi)|\left(\frac{1-r}{1-\rho}\right)^{c_{1}-1}
\left(\frac{1+\rho}{1+r}\right)^{c_{1}+1}.
\eeq
Applying (\ref{eqt-10}) and the arbitrariness of $\mu$, we see that
\beq\label{eqt-11}
\Lambda_{f}(r\xi)&\geq&\Lambda_{f}(\rho\xi)\left(\frac{1-r}{1-\rho}\right)^{c_{1}-1}
\left(\frac{1+\rho}{1+r}\right)^{c_{1}+1}\\ \nonumber
&\geq&\frac{1}{2^{1+c_{1}}}\Lambda_{f}(\rho\xi)\left(\frac{1-r}{1-\rho}\right)^{c_{1}-1}.
\eeq

Now, we  prove that $|b_{2}|\leq\frac{c}{2}.$ Since
$\omega(0)=0$ and  $|\omega(z)|\leq c$ in $\ID$, follows that
$|\omega'(0)|\leq c$ and hence, from
$$\omega(z)=\frac{g'(z)}{h'(z)}=\frac{2b_{2}z+3b_{3}z^{2}+\cdots}{1+2a_{2}z+\cdots}
=2b_{2}z+(3b_{3}-4a_{2}b_{2})z^{2}+\cdots,
$$
we obtain that $|\omega'(0)|=|2b_{2}|\leq c.$ Then
\be\label{eq-20}
|b_{2}|\leq\frac{c}{2}.
\ee

Finally, we prove the sharpness part.  For $z\in\mathbb{D}$, let
$f(z)=h(z)+\overline{g(z)},$ where $h(z)=z$ and
$g(z)=\frac{c}{2}\overline{z}^{2}$. For $w\in h(\mathbb{D})=\ID,$ let
$$F(w)=f(h^{-1}(w))=w+\overline{g(h^{-1}(w))}=w+\frac{c}{2}\overline{w}^2.
$$
Since $h(\mathbb{D})$ is convex, we see that for any $w_{1},
w_{2}\in h(\mathbb{D})$,
\beq\label{eq-y1}
l(F(\gamma))&\leq&\int_{\gamma}\left|F_{w}(w)\,dw+F_{\overline{w}}(w)\,d\overline{w}\right|\\
\nonumber
&\leq&\int_{\gamma}(|F_{w}(w)|+|F_{\overline{w}}(w)|)|dw|\\
\nonumber &\leq&(1+c)|w_{2}-w_{1}|,
\eeq
where $\gamma$ is a line segment joining $w_{1}$ and $w_{2}$. On the other hand,
\beq\label{eq-y2}
|F(w_{2})-F(w_{1})|&\geq&|w_{2}-w_{1}|-\frac{c}{2}|w^{2}_{2}-w_{1}^{2}|\\
\nonumber &\geq&(1-c)|w_{2}-w_{1}|.
\eeq
Equations (\ref{eq-y1}) and (\ref{eq-y2}) yield that
$$l(F(\gamma))\leq\frac{1+c}{1-c}|F(w_{2})-F(w_{1})|.
$$
Hence $f(\mathbb{D})$ is $\left(\frac{1+c}{1-c}\right)$-linearly
connected, where $f(z)=z+\frac{c}{2}\overline{z}^{2}$.  Therefore,
the extremal function $f(z)=z+\frac{c}{2}\overline{z}^{2}$ shows
that the estimate of (\ref{eq-20}) is sharp.

Next we prove the first part of  (III).
 Let
$T_{\theta}=h+e^{i\theta}g$, where $\theta\in[0,2\pi)$. First of
all, we prove that $T_{\theta}$ is univalent and
$T_{\theta}(\mathbb{D})$ is a $\frac{M(1+c)}{1-c(1+2M)}$-linearly
connected domain. For $w\in f(\mathbb{D})$, let
$$H(w)=T_{\theta}(f^{-1}(w))=w-\overline{G(w)}+e^{i\theta}G(w),
$$
where $G(w)=g(f^{-1}(w))$. By the chain rule, we get
$$G_{w}=g'\cdot(f^{-1})_{w}~\mbox{ and }~G_{\overline{w}}=g'\cdot(f^{-1})_{\overline{w}}.
$$
Differentiating both sides of equation $f^{-1}(f(z))=z$ yields the
relations
$$(f^{-1})_{w}\cdot h'+(f^{-1})_{\overline{w}}\cdot g'=1~\mbox{ and }~(f^{-1})_{w}\cdot
\overline{g'}+(f^{-1})_{\overline{w}}\cdot \overline{h'}=0,
$$
which imply that
$$(f^{-1})_{w}=\frac{\overline{h'}}{J_{f}}~\mbox{ and }~(f^{-1})_{\overline{w}}=-\frac{\overline{g'}}{J_{f}}.
$$
Since $|\omega|\leq c$ (by hypotheses), we see that
\be\label{eq-14}
\Lambda_{G}(w)=\frac{|\omega|}{1-|\omega|}\leq\frac{c}{1-c}.
\ee
Since $f(\mathbb{D})$ is $M$-linearly connected, we know that for
any two points $w_{1},w_{2}\in f(\mathbb{D})$, there is  a curve
$\gamma\subset f(\mathbb{D})$ joining $w_{1}$ and $w_{2}$ such that
$l(\gamma)\leq M|w_{1}-w_{2}|$. Now, we set $\Gamma=H(\gamma)$. By
elementary calculations, we have
\be\label{eq-15}
H_{w}=1-\overline{G_{\overline{w}}}+e^{i\theta}G_{w}
~\mbox{ and }~H_{\overline{w}}=-\overline{G_{w}}+e^{i\theta}G_{\overline{w}}.
\ee
By using (\ref{eq-14}) and (\ref{eq-15}), we get
\beq\label{eq-16}
l(\Gamma)&=&\int_{\gamma}\left|H_{w}(w)\,dw+H_{\overline{w}}(w)\,d\overline{w}\right|\\
\nonumber &\leq&\int_{\gamma}\left(1+2\Lambda_{G}(w)\right)|dw|\\
\nonumber &\leq&\left(1+\frac{2c}{1-c}\right)l(\gamma)\\
\nonumber&\leq&M\left(\frac{1+c}{1-c}\right)|w_{1}-w_{2}|.
\eeq
On the other hand, we have
\be\label{eq-17}
|H(w_{2})-H(w_{1})|\geq|w_{2}-w_{1}|-2\int_{\gamma}\Lambda_{G}(w)|dw|\geq|w_{1}-w_{2}|\left(1-\frac{2cM}{1-c}\right),
\ee
which shows that for all $\theta\in[0,2\pi)$, $T_{\theta}$ is univalent.
Equations (\ref{eq-16}) and (\ref{eq-17}) yield that
$$l(\Gamma)\leq\frac{M(1+c)}{1-c(1+2M)}|H(w_{2})-H(w_{1})|,
$$
which implies that $T_{\theta}(\mathbb{D})$ is a $\frac{M(1+c)}{1-c(1+2M)}$-linearly
connected domain.

By \cite[Proposition 5.6]{P}, we know that $T_{\theta}$ is
continuous in $\overline{\mathbb{D}}$ with values in
$\mathbb{C}\cup\{\infty\}.$ Applying \cite[Theorem 5.7 (5)]{P} to
$T_{\theta}$, we see that there is a positive constant $c_{2}$ and $c_{3}<2$ such that
for $\zeta_{1}, \zeta_{2}\in\partial\mathbb{D}$,
\be\label{eq-18}
|T_{\theta}(\zeta_{1})-T_{\theta}(\zeta_{2})|\geq c_{2}|\zeta_{1}-\zeta_{2}|^{c_{3}}.
\ee
Inequality (\ref{eq-18}) and the arbitrariness of $\theta\in[0,2\pi)$ gives
$$|f(\zeta_{1})-f(\zeta_{2})|\geq c_{2}|\zeta_{1}-\zeta_{2}|^{c_{3}}.
$$
This completes the proof of the first part of  (III).

Next we show the second part  of (III). By the proof of the first
part  in (III), we know that for all $\theta\in[0,2\pi)$,
$T_{\theta}=h+e^{i\theta}g$ is univalent in $\mathbb{D}$. Applying
the result of  de Branges in \cite{B}, we see that for all
$\theta\in[0,2\pi)$,
$$|a_{n}+e^{i\theta}b_{n}|\leq n,
$$
which implies that $|a_{n}|+|b_{n}|\leq n,$ for $n\geq2.$
\hfill $\Box$ 


\subsection*{Proof of Theorem \ref{thm-1}}
 We first prove (I). Without loss of generality, we assume that
$F=h-g$ is $\alpha$-close-to convex. By \cite[Proposition 5.8]{P},
we know that $\Omega=F(\mathbb{D})$ is $M^{\ast}$-linearly connected
domain, where $M^{\ast}=\frac{1}{\cos\frac{\alpha\pi}{2}}$. Let
$$H(w)=h(F^{-1}(w))=w+g(F^{-1}(w)).
$$
For any two distinct $w_{1},w_{2}\in\Omega$, by hypothesis, there
is a curve $\gamma\subset\Omega$ joining $w_{1}$ and $w_{2}$ such
that $l(\gamma)\leq M^{\ast}|w_{1}-w_{2}|$. Set $\Gamma=H(\gamma).$
Then we have
\beq\label{eqt-2}
l(\Gamma)&=&\int_{\Gamma}|dH(w)|=\int_{\Gamma}\left|H'(w)\,dw\right|\\
\nonumber
&=&\int_{\gamma}\left|1+\frac{g'(F^{-1}(w))}{h'(F^{-1}(w))-g'(F^{-1}(w))}\right||dw|\\
\nonumber
&\leq&\int_{\gamma}\left(1+\frac{|\omega(F^{-1}(w))|}{1-|\omega(F^{-1}(w))|}\right)|dw|\\
\nonumber &\leq&\frac{l(\gamma)}{1-M_{1}}\\ \nonumber &\leq&
\frac{M^{\ast}}{1-M_{1}}|w_{1}-w_{2}|.
\eeq
On the other hand,
\beq\label{eqt-3}
|H(w_{2})-H(w_{1})|&\geq&|w_{2}-w_{1}|-\left|g(F^{-1}(w_{2}))-g(F^{-1}(w_{1}))\right|\\
\nonumber
&\geq&|w_{2}-w_{1}|-\int_{\gamma}\left|\big(g(F^{-1}(w))\big)'\right||dw|\\
\nonumber&=&\frac{1-M_{1}(1+M^{\ast})}{1-M_{1}}|w_{2}-w_{1}|.
\eeq
By (\ref{eqt-2}) and (\ref{eqt-3}), we obtain
$$l(\Gamma)\leq \frac{M^{\ast}}{1-M_{1}(1+M^{\ast})}|H(w_{2})-H(w_{1})|,
$$
which shows that $h(\mathbb{D})$ is $M_{2}$-linearly connected domain,
where
\[
M_{2}=\frac{1}{\cos\frac{\alpha\pi}{2}-M_{1}(1+\cos\frac{\alpha\pi}{2})}.
\]
The univalency of $h$ follows from (\ref{eqt-3}).

Next we  prove (II). Define
$$T(w)=f_{\theta}(F^{-1}(w))=w+g(F^{-1}(w))+e^{i\theta}\overline{g(F^{-1}(w))}.
$$
Let $\Gamma_{1}=T(\gamma)$. Then we find that
\beq\label{eqt-4}
l(\Gamma)&=&\int_{\Gamma}|dT(w)|=\int_{\Gamma}\left|T_{w}(w)\,dw+T_{\overline{w}}(w)\,d\overline{w}\right|\\
\nonumber
&\leq&\int_{\gamma}\left(1+\frac{2|g'(F^{-1}(w))|}{|h'(F^{-1}(w))-g'(F^{-1}(w))|}\right)|dw|
\\ \nonumber
&\leq&\int_{\gamma}\left(1+\frac{2|\omega(F^{-1}(w))|}{1-|\omega(F^{-1}(w))|}\right)|dw|
\\ \nonumber
&\leq&\frac{1+M_{3}}{1-M_{3}}l(\gamma) \leq
\frac{M^{\ast}(1+M_{3})}{1-M_{3}}|w_{1}-w_{2}|
\eeq
and
\beq\label{eqt-5}
|T(w_{2})-T(w_{1})|&\geq&|w_{2}-w_{1}|-2\left|g(F^{-1}(w_{2}))-g(F^{-1}(w_{1}))\right|\\
\nonumber
&\geq&|w_{2}-w_{1}|-2\int_{\gamma}\left|\big(g(F^{-1}(w))\big)'\right||dw|\\
\nonumber&=&\frac{1-M_{3}(1+2M^{\ast})}{1-M_{3}}|w_{2}-w_{1}|.
\eeq
By (\ref{eqt-4}) and (\ref{eqt-5}), we obtain
$$l(\Gamma_{1})\leq \frac{M^{\ast}(1+M_{3})}{1-M_{3}(1+2M^{\ast})}|T(w_{2})-T(w_{1})|,
$$
which implies that $f_{\theta}(\mathbb{D})$ is $M_{4}$-linearly
connected domain, where $\theta\in[0,2\pi)$ and
\[
M_{4}=\frac{1+M_{3}}{\cos\frac{\alpha\pi}{2}-M_{3}(2+\cos\frac{\alpha\pi}{2})}.
\]
The univalency of $f_{\theta}$ follows from (\ref{eqt-5}). Since for
$z\in\mathbb{D}$,
$$\frac{\Lambda_{f_{\theta}}(z)}{\lambda_{f_{\theta}}(z)}\leq\frac{1+M_{3}}{1-M_{3}},
$$
we see that $f_{\theta}$ is a $K$-quasiconformal harmonic mapping, where $K=\frac{1+M_{3}}{1-M_{3}}$.
\hfill $\Box$ 



\normalsize

\end{document}